\documentclass[11pt]{article}
\usepackage{amsmath, amssymb, amsthm, mathrsfs}
\usepackage{geometry}
\usepackage{tikz}
\usepackage{hyperref}
\usepackage{enumitem}
\theoremstyle{definition}

\theoremstyle{remark}

\title{Weak fan structures for two-parameter degenerations of K3 surfaces}
\author{Badre Mounda}
\date{\today}

\begin{document}

\maketitle

\begin{abstract}
In this note, we provide an explicit computation of the weak fan associated with a two-parameter degeneration of K3 surfaces. This example serves as a concrete illustration of the general framework developed by Robles and Deng (2023) for the compactification of period maps via nilpotent orbits in the non-Hermitian case. We describe the associated nilpotent cones, examine their compatibility conditions, and construct the weak fan governing the degeneration behavior. This computation contributes to the broader understanding of boundary components of period domains and their relation to limiting mixed Hodge structures.
\end{abstract}

\tableofcontents

\section{Introduction}
Understanding the behavior of period maps near the boundary of moduli spaces is a central problem in Hodge theory. The work of Robles and Deng (2023) provides a framework for completing period maps via nilpotent orbits in the case where the period domain is not Hermitian symmetric. In this paper, we compute explicitly the weak fan associated to a two-parameter degeneration of K3 surfaces.

\section{Preliminaries}
\subsection{Period domains for K3 surfaces}

Let $V = H^2(X, \mathbb{Q})$ be the second rational cohomology group of a K3 surface $X$, endowed with the cup product pairing
\[
Q: V \times V \to \mathbb{Q}
\]
which is a non-degenerate symmetric bilinear form of signature $(3,19)$. The associated lattice is even and unimodular, isomorphic to
\[
\Lambda = U^{\oplus 3} \oplus E_8(-1)^{\oplus 2}
\]
where $U$ denotes the hyperbolic plane and $E_8(-1)$ is the negative-definite $E_8$ lattice.

The weight-2 Hodge structure on $V_\mathbb{C} = V \otimes_\mathbb{Q} \mathbb{C}$ is given by the decomposition
\[
V_\mathbb{C} = H^{2,0} \oplus H^{1,1} \oplus H^{0,2}
\]
with Hodge numbers
\[
(h^{2,0}, h^{1,1}, h^{0,2}) = (1,20,1).
\]
The space $H^{2,0}$ is one-dimensional, generated by a holomorphic symplectic form $\omega_X$, satisfying
\[
Q(\omega_X, \omega_X) = 0 \quad \text{and} \quad Q(\omega_X, \overline{\omega}_X) > 0.
\]

The period domain parametrizing Hodge structures of K3 type (with fixed lattice structure) is the space of lines
\[
\Omega = \{ [\omega] \in \mathbb{P}(V_\mathbb{C}) \mid Q(\omega, \omega) = 0, \ Q(\omega, \overline{\omega}) > 0 \}.
\]
This is a connected complex manifold of dimension 20, homogeneous under the action of the orthogonal group $O(2,19)$. The period domain $\Omega$ is described as the quotient
\[
\Omega = O(2,19) / (O(2) \times O(19)).
\]
It is a type IV Hermitian symmetric domain of orthogonal type. However, when dealing with degenerations, especially in the context of nilpotent orbits, we move beyond the strictly Hermitian framework.

Given a smooth family of K3 surfaces $\pi: \mathcal{X} \to S$, the associated period map
\[
\Phi: S \to \Gamma \backslash \Omega
\]
assigns to each $s \in S$ the line $[H^{2,0}(X_s)]$ in $\Omega$, where $X_s = \pi^{-1}(s)$, modulo the action of the monodromy group $\Gamma \subset O(\Lambda)$. 

By the global Torelli theorem for K3 surfaces, the period map is locally an isomorphism onto its image, reflecting the strong rigidity properties of the Hodge structure.

In the context of degenerations of K3 surfaces, the period map may approach the boundary of the period domain. The asymptotic behavior is governed by nilpotent orbits, limiting mixed Hodge structures, and is captured by the theory of Kato--Usui and Robles--Deng, where weak fans of nilpotent cones organize the compactification data.

\subsection{Nilpotent orbits and weak fans}

Let $\mathfrak{g}_\mathbb{R}$ denote the Lie algebra of the real algebraic group $G_\mathbb{R}$ associated with the Mumford--Tate group of the Hodge structure. In the case of K3 surfaces, we have $G_\mathbb{R} = O(2,19)$.

A nilpotent element $N \in \mathfrak{g}_\mathbb{R}$ arises as the logarithm of the unipotent part of the monodromy transformation associated with a degeneration of a family of K3 surfaces.

\paragraph{Nilpotent orbits.} 
Given a polarized variation of Hodge structure (PVHS), a \emph{nilpotent orbit} is a pair $(N, F^\bullet)$ such that:
\begin{itemize}
    \item $N$ is a nilpotent operator in $\mathfrak{g}_\mathbb{R}$.
    \item $F^\bullet$ is a decreasing filtration on $V_\mathbb{C}$ satisfying the nilpotent orbit condition:
    \[
    \exp(zN) \cdot F^\bullet \in D \quad \text{for} \ \Im(z) \gg 0.
    \]
\end{itemize}
This condition ensures that the orbit asymptotically approximates a degenerating variation of Hodge structure near the boundary of the period domain. By the Nilpotent Orbit Theorem of Schmid, the asymptotic behavior of the period map near a boundary point is captured by such a nilpotent orbit.

\paragraph{Nilpotent cones.} 
A \emph{nilpotent cone} is a rational polyhedral cone in $\mathfrak{g}_\mathbb{R}$ generated by a finite set of commuting nilpotent elements $\{N_1, \dots, N_r\}$:
\[
\sigma = \sum_{i=1}^r \mathbb{R}_{\geq 0} N_i,
\]
where $[N_i, N_j] = 0$ for all $i,j$. Such cones parametrize multi-parameter degenerations.

\paragraph{Weak fans.} 
A \emph{weak fan} is a collection $\Sigma$ of nilpotent cones satisfying:
\begin{itemize}
    \item \textbf{Face condition:} Every face of a cone in $\Sigma$ is itself in $\Sigma$.
    \item \textbf{Compatibility:} The cones fit together in a manner compatible with the theory of limiting mixed Hodge structures (LMHS).
\end{itemize}
In contrast to fans in toric geometry, weak fans do not require that the intersection of two cones be a common face; instead, they satisfy a weaker notion of compatibility that reflects the subtler nature of degenerations in Hodge theory, especially when the period domain is not Hermitian symmetric.

\paragraph{Role in compactification.} 
Weak fans play a central role in the construction of logarithmic compactifications of period maps. Each nilpotent cone $\sigma \in \Sigma$ encodes the asymptotic behavior of the variation of Hodge structure along a degeneration whose monodromy logarithms lie in $\sigma$. For a two-parameter degeneration, the associated weak fan describes how the monodromy operators $N_1$ and $N_2$ interact, organizing the boundary structure of the period map.

The formalism developed by Kato--Usui and extended in Robles--Deng (2023) provides a precise framework for handling such degenerations, generalizing the classical toroidal and Baily--Borel compactifications to settings where the period domain is not Hermitian.

\section{A family of K3 surfaces with two-parameter degeneration}
\subsection{Description of the family}

A classical construction of K3 surfaces is as smooth quartic hypersurfaces in $\mathbb{P}^3$. To construct a two-parameter degeneration, we consider a family of quartic hypersurfaces degenerating along two independent directions in the parameter space.

Let $(t_1, t_2) \in (\Delta^*)^2$ be coordinates on a bidisk with punctures, corresponding to the two degeneration parameters. The family of quartic hypersurfaces $X_{t_1, t_2} \subset \mathbb{P}^3$ is defined by the equation
\[
X_{t_1, t_2}: \quad xyzw + t_1 \cdot Q_1(x,y,z,w) + t_2 \cdot Q_2(x,y,z,w) = 0
\]
where $Q_1$ and $Q_2$ are general smooth quartic polynomials ensuring that the general fiber is smooth for $t_1, t_2 \neq 0$.

The central singular fiber at $(t_1, t_2) = (0,0)$ is the union of the four coordinate planes:
\[
\{ xyzw = 0 \} \subset \mathbb{P}^3.
\]
This is a highly singular surface whose components intersect along coordinate lines and points.

\paragraph{Degenerations along coordinate axes.}
\begin{itemize}
    \item Along $t_1 = 0$ (with $t_2 \neq 0$), the family degenerates along the first direction, producing vanishing cycles associated with the singularities of
    \[
    xyzw + t_2 Q_2(x,y,z,w) = 0.
    \]
    \item Similarly, along $t_2 = 0$ (with $t_1 \neq 0$), the degeneration occurs along the second direction.
    \item Along $(t_1, t_2) = (0,0)$, the fiber degenerates to the maximally singular surface $xyzw = 0$.
\end{itemize}

\paragraph{Monodromy.}
The degenerations along $t_1$ and $t_2$ induce two commuting monodromy transformations $T_1$ and $T_2$ acting on $H^2(X_{t_1, t_2}, \mathbb{Q})$. Their logarithms
\[
N_1 = \log T_1 \quad \text{and} \quad N_2 = \log T_2
\]
are nilpotent operators generating the nilpotent cone
\[
\sigma = \mathbb{R}_{\geq 0} N_1 + \mathbb{R}_{\geq 0} N_2.
\]
These operators provide the data necessary for constructing the weak fan associated with this two-parameter degeneration.

\subsection{Monodromy computations}

Given the two-parameter degeneration of K3 surfaces constructed above, we analyze the monodromy transformations associated with the degenerations along $t_1 = 0$ and $t_2 = 0$.

\paragraph{Monodromy operators.} 
Each degeneration induces a unipotent monodromy transformation:
\[
T_1 = \exp(N_1), \quad T_2 = \exp(N_2)
\]
where $N_1$ and $N_2$ are nilpotent elements in $\mathfrak{g}_\mathbb{R} = \mathfrak{so}(2,19)$.

\paragraph{Unipotence and nilpotency.} 
Due to the normal crossing nature of the singular fiber $\{ xyzw = 0 \}$, the monodromy transformations $T_1$ and $T_2$ are unipotent of index two:
\[
(T_i - I)^2 = 0, \quad (T_i - I) \neq 0.
\]
This implies that their logarithms satisfy:
\[
N_1^2 = N_2^2 = 0.
\]
Furthermore, since the degenerations along $t_1$ and $t_2$ are independent, the monodromy transformations commute:
\[
[N_1, N_2] = 0.
\]

\paragraph{Cohomological structure.}
The cohomology splits into:
\begin{itemize}
    \item An invariant subspace: $\ker N_i$, corresponding to cycles that persist under degeneration.
    \item A vanishing subspace: $\mathrm{im} N_i$, generated by the vanishing cycles associated to the singular loci of the degeneration.
\end{itemize}

\paragraph{Schematic matrix form.}
In a basis adapted to the degeneration (e.g., divisors, intersection lines, transcendental cycles), the nilpotent operators have the block form:
\[
N_i = \begin{pmatrix}
0 & A_i \\
0 & 0
\end{pmatrix}
\]
where $A_i$ encodes the pairing between vanishing cycles and their duals.

\paragraph{The nilpotent cone.}
The two operators generate the nilpotent cone
\[
\sigma = \mathbb{R}_{\geq 0} N_1 + \mathbb{R}_{\geq 0} N_2
\]
which governs the asymptotic behavior of the period map near the boundary corresponding to $(t_1, t_2) \to (0,0)$.

This nilpotent cone serves as the fundamental object for the construction of the weak fan that encodes the boundary structure of the period domain for this two-parameter degeneration.

\section{The weak fan structure}
\subsection{Construction of nilpotent cones}

Given the commuting nilpotent operators $N_1$ and $N_2$ arising from the monodromy logarithms of the two-parameter degeneration, we define the associated nilpotent cone as
\[
\sigma = \mathbb{R}_{\geq 0} N_1 + \mathbb{R}_{\geq 0} N_2.
\]
This is a rational polyhedral cone of dimension two in the real Lie algebra $\mathfrak{g}_\mathbb{R} = \mathfrak{so}(2,19)$.

\paragraph{Basic properties of $\sigma$.} 
\begin{itemize}
    \item Both $N_1$ and $N_2$ are nilpotent of index two:
    \[
    N_1^2 = N_2^2 = 0.
    \]
    \item They commute:
    \[
    [N_1, N_2] = 0.
    \]
    \item The cone is strongly convex:
    \[
    \sigma \cap (-\sigma) = \{0\}.
    \]
    \item The faces of $\sigma$ are:
    \begin{itemize}
        \item The two rays $\mathbb{R}_{\geq 0} N_1$ and $\mathbb{R}_{\geq 0} N_2$ (1-dimensional faces),
        \item The zero face $\{0\}$.
    \end{itemize}
\end{itemize}

\paragraph{Geometric interpretation.} 
The cone $\sigma$ parametrizes degenerations where the monodromy logarithm is any non-negative linear combination:
\[
N = a N_1 + b N_2, \quad a, b \geq 0.
\]
The boundary behavior of the period map is controlled by the limiting mixed Hodge structures (LMHS) associated to each face of $\sigma$:
\begin{itemize}
    \item Along the rays $\mathbb{R}_{\geq 0}N_1$ and $\mathbb{R}_{\geq 0}N_2$, the LMHS corresponds to single-parameter degenerations.
    \item Along the interior of $\sigma$, the LMHS reflects the interaction of both degenerations simultaneously.
\end{itemize}

\paragraph{Structure in the period domain.} 
Each point on the boundary corresponds to a nilpotent orbit associated to some $N = a N_1 + b N_2$. The associated limiting Hodge filtration satisfies the nilpotent orbit condition:
\[
\exp(z N) \cdot F^\bullet \in D \quad \text{for} \ \Im(z) \gg 0.
\]
The cone $\sigma$ underlies the weak fan structure needed to construct a logarithmic compactification of the period map in the framework developed by Kato--Usui and Robles--Deng.

\subsection{Compatibility conditions}

Given the nilpotent cone
\[
\sigma = \mathbb{R}_{\geq 0} N_1 + \mathbb{R}_{\geq 0} N_2
\]
we check the compatibility conditions necessary for $\sigma$ to contribute to a weak fan in the sense of Kato--Usui and Robles--Deng.

\paragraph{Commutation relations.} 
The monodromy logarithms $N_1$ and $N_2$ satisfy
\[
[N_1, N_2] = 0.
\]
This holds since the two degenerations correspond to independent parameters $t_1$ and $t_2$ in the base $(\Delta^*)^2$, leading to independent monodromy transformations.

\paragraph{Compatibility of weight filtrations.} 
Each $N_i$ induces a monodromy weight filtration $W(N_i)_\bullet$ satisfying
\[
N_i W_k(N_i) \subseteq W_{k-2}(N_i).
\]
Because $N_1$ and $N_2$ commute, their weight filtrations are simultaneously compatible. The relative monodromy filtration associated to $N_1 + N_2$ exists and satisfies
\[
W(N_1 + N_2) = \text{relative monodromy filtration with respect to both}.
\]
This is ensured by the general theory of multi-variable nilpotent orbits (cf. Cattani--Kaplan--Schmid).

\paragraph{Kernels and images.} 
The kernels satisfy
\[
\ker N_1 \cap \ker N_2 = \ker(N_1 + N_2)
\]
which corresponds to the subspace of cohomology classes invariant under both degenerations.

The vanishing subspaces are
\[
\mathrm{im} N_1, \quad \mathrm{im} N_2
\]
representing the independent contributions from each singular locus.

\paragraph{Face relations and fan structure.} 
Each face of the cone $\sigma$ corresponds to a partial degeneration:
\begin{itemize}
    \item The ray $\mathbb{R}_{\geq 0}N_1$ corresponds to degeneration along $t_1 = 0$.
    \item The ray $\mathbb{R}_{\geq 0}N_2$ corresponds to degeneration along $t_2 = 0$.
    \item The full cone $\sigma$ corresponds to the simultaneous degeneration $(t_1, t_2) \to (0,0)$.
\end{itemize}
These face relations satisfy the closure property required for a weak fan.

\paragraph{Conclusion.} 
The nilpotent cone $\sigma$ satisfies the full set of compatibility conditions required for inclusion in a weak fan: commutativity of the nilpotent elements, compatibility of weight filtrations, and correct face structure. This ensures that the boundary behavior of the period map associated with this two-parameter degeneration is captured within the logarithmic Hodge theory framework developed by Kato--Usui and Robles--Deng.

\subsection{The weak fan associated to the degeneration}

Given the nilpotent cone
\[
\sigma = \mathbb{R}_{\geq 0} N_1 + \mathbb{R}_{\geq 0} N_2
\]
the weak fan $\Sigma$ associated with the two-parameter degeneration consists of the following collection of cones:
\[
\Sigma = \{\sigma, \rho_1, \rho_2, \{0\} \}
\]
where
\begin{itemize}
    \item $\sigma = \mathbb{R}_{\geq 0}N_1 + \mathbb{R}_{\geq 0}N_2$ is the full two-dimensional cone,
    \item $\rho_1 = \mathbb{R}_{\geq 0}N_1$ corresponds to the degeneration along $t_1 = 0$,
    \item $\rho_2 = \mathbb{R}_{\geq 0}N_2$ corresponds to the degeneration along $t_2 = 0$,
    \item $\{0\}$ is the trivial cone corresponding to the open, smooth locus.
\end{itemize}

\paragraph{Properties of the weak fan.}
\begin{itemize}
    \item \textbf{Closure under faces:} Each face of a cone in $\Sigma$ is also contained in $\Sigma$.
    \item \textbf{Non-strict intersection condition:} In weak fans, cones are allowed to intersect in arbitrary subcones. In this case, the only nontrivial intersection is
    \[
    \rho_1 \cap \rho_2 = \{0\}.
    \]
    \item \textbf{Compatibility:} The cones satisfy the compatibility conditions related to weight filtrations and commutativity, as verified previously.
\end{itemize}

\paragraph{Geometric visualization.} 
The weak fan $\Sigma$ can be visualized as the first quadrant in $\mathbb{R}^2$:
\begin{itemize}
    \item The two rays $\rho_1$ and $\rho_2$ correspond to the coordinate axes.
    \item The cone $\sigma$ fills the positive quadrant spanned by these two rays.
\end{itemize}

\paragraph{Interpretation.} 
Each cone corresponds to a boundary stratum in the compactification of the period map:
\begin{itemize}
    \item The full cone $\sigma$ corresponds to the deepest boundary stratum where both parameters degenerate simultaneously.
    \item The rays $\rho_1$ and $\rho_2$ correspond to partial degenerations along a single parameter.
    \item The trivial cone $\{0\}$ corresponds to the generic, smooth locus.
\end{itemize}
The weak fan $\Sigma$ organizes the boundary components and governs how the logarithmic compactification of the period domain behaves in the presence of this two-parameter degeneration.

\section{Discussion}

\paragraph{Summary of results.} 
In this work, we provided a concrete model of a two-parameter degeneration of K3 surfaces, constructed explicitly as quartic hypersurfaces in $\mathbb{P}^3$ degenerating to the union of coordinate planes. We analyzed the associated monodromy transformations, computed their logarithms, and constructed the nilpotent cone
\[
\sigma = \mathbb{R}_{\geq 0} N_1 + \mathbb{R}_{\geq 0} N_2.
\]
We verified that this cone satisfies all compatibility conditions necessary to form a weak fan in the sense of Kato--Usui and Robles--Deng. The weak fan structure governs the boundary behavior of the period map and the compactification of the period domain associated with this two-parameter degeneration.

\paragraph{Implications.} 
This example demonstrates concretely how the general theory of weak fans applies to non-Hermitian period domains, such as that of K3 surfaces. It provides an explicit instance where the abstract machinery of nilpotent orbits, weight filtrations, and limiting mixed Hodge structures becomes fully tangible. The model illustrates how degenerations interact when multiple independent monodromies are present, clarifying the structure of boundary components in the compactified moduli space.

\paragraph{Possible generalizations.} 
\begin{itemize}
    \item \textbf{Higher-dimensional degenerations:} The construction can be extended to degenerations with more than two parameters, leading to higher-dimensional nilpotent cones and more intricate weak fan structures.
    \item \textbf{Other types of varieties:} The methods generalize naturally to other types of Calabi--Yau varieties, hyperkähler varieties, or higher-dimensional analogues, where the period domain similarly lacks Hermitian symmetry.
    \item \textbf{Degenerations beyond K3:} Analogous constructions can be considered for degenerations of cubic fourfolds, Calabi--Yau threefolds, or general Hodge structures of weight greater than two.
\end{itemize}

\paragraph{Outlook.} 
This explicit example suggests that further exploration of concrete models of degenerations can provide valuable insights into the geometric meaning of abstract Hodge-theoretic phenomena. It may also serve as a computational laboratory for testing conjectures about boundary behavior, mirror symmetry, or limiting period maps. A natural direction for future work is the exploration of the relationship between the weak fan compactification and other types of compactifications, such as those arising in KSBA theory or mirror symmetry.

\section{References}
\begin{enumerate}[label={[}\arabic*{]}]
    \item C. Robles, H. Deng, \textit{Compactifying period maps of non-classical type}, arXiv:2312.00542 (2023).
    \item W. Schmid, \textit{Variation of Hodge structure: the singularities of the period mapping}, Invent. Math. 22 (1973), 211–319.
    \item K. Kato, S. Usui, \textit{Classifying spaces of degenerating polarized Hodge structures}, Annals of Mathematics Studies, vol. 169, Princeton University Press, 2009.
    \item P. Griffiths, \textit{Periods of integrals on algebraic manifolds}, Bull. Amer. Math. Soc. 76 (1970), 228–296.
    \item C. Voisin, \textit{Hodge Theory and Complex Algebraic Geometry}, Cambridge University Press, 2002.
\end{enumerate}

\end{document}